 \documentstyle[graphicx,amscd,amssymb,verbatim,12pt,righttag,color]{amsart}
 \newtheorem{theorem}{Theorem}[section]
 
 \newtheorem{Prop}[theorem]{Proposition}
 \newtheorem{Lem}[theorem]{Lemma}
 \newtheorem{Cor}[theorem]{Corollary}
 \newtheorem{Rem}[theorem]{Remark}
 \newtheorem{Example}[theorem]{Example}

\newcommand{\eproof}{\hfill$\square$}

 \numberwithin{equation}{section}
 \renewcommand{\rm}{\normalshape}
 
\begin{document}

\title  {Spectral measures associated with the factorization of the Lebesgue measure on a set via convolution}

\author{Jean-Pierre Gabardo}
\email{gabardo@@mcmaster.ca}

\author{Chun-Kit Lai}
 \email{cklai@@math.mcmaster.ca}

\address{Department of Mathematics and Statistics, McMaster University,
Hamilton, Ontario, L8S 4K1, Canada}

\thanks{The first named author was supported by an NSERC grant.}

\begin{abstract}
Let $Q$ be a fundamental domain of some full-rank lattice in ${\Bbb R}^d$ and let $\mu$ and $\nu$ be two positive Borel measures on ${\Bbb R}^d$ such that the convolution $\mu\ast\nu$ is a multiple of $\chi_Q$. We consider the problem as to whether or not both measures must be spectral (i.e. each of their respective associated $L^2$ space admits an orthogonal basis of exponentials) and we show that this is the case when $Q = [0,1]^d$. This theorem yields a large class of examples of spectral measures which are either absolutely continuous, singularly continuous or purely discrete spectral measures. In addition, we propose a
generalized Fuglede's conjecture for spectral measures on ${\Bbb R}^1$ and we show that it implies the classical
Fuglede's conjecture on ${\Bbb R}^1$.
\end{abstract}
\maketitle 

\section{Introduction}

Let $\mu$ be a compactly supported Borel probability measure on
${\Bbb R}^d$. We say that $\mu$ is a {\it spectral measure} if there
exists a countable set $\Lambda\subset {\Bbb R}^d$ called {\it spectrum} such  that
$E(\Lambda): = \{e^{2\pi i \langle\lambda,x\rangle}:
\lambda\in\Lambda\}$ is an orthonormal basis for $L^2(\mu)$. If $\Omega\subset{\Bbb R}^d$ is measurable with finite positive Lebesgue measure and
$d\mu(x) =\chi_{\Omega}(x)dx$ is a spectral measure, then we say that $\Omega$
is a {\it spectral set}. Spectral sets were first introduced by Fuglede (\cite{[Fu]}) and
 have a very delicate and mysterious relationship with translational tiling because of the
 {\it spectral set conjecture} (known also as {\it Fuglede's conjecture}) proposed by Fuglede.

 \medskip

 \noindent{\bf Conjecture (Fuglede's Conjecture):} {\it A bounded measurable set $\Omega$
on ${\Bbb R}^d$ of positive Lebesgue measure is a spectral set if and only if
$\Omega$ is a translational tile.}

 \medskip

 We say that $\Omega$ is a translational tile if there exists a
discrete set ${\mathcal J}$ such that $\bigcup_{t\in{\mathcal J}}(\Omega+t) = {\Bbb R}^d$,
and the Lebesgue measure of  $(\Omega+t)\cap(\Omega+t')$ is zero for any distinct
 $t$ and $t'$ in ${\mathcal J}$.  Although this conjecture was eventually disproved
in dimension $d\geq 3$ (\cite{[T], [KM], [KM1]}),
most of the known examples of spectral sets are constructed from
translational tiles. An important class of examples of spectral sets
constructed in \cite{[PW]} consists of sets of the form $A+[0,1]$  tiling $[0,N]$ for some
$N$, where $A\subset\Bbb Z$. In fact, in this case, the corresponding
equally weighted discrete measure on $A$ is a spectral measure.

\medskip

The first singular spectral measure was constructed by Jorgensen and Pedersen \cite{[JP]}.
They showed that the standard
 Cantor measures are spectral measures if the contraction is $\frac{1}{2n}$, while there are at most two orthogonal
 exponentials when the contraction
 is $\frac{1}{2n+1}$. Following this discovery, more spectral self-similar/self-affine measures were also found (\cite{[S]}, \cite{[LaW]}, \cite{[DJ]}). In these investigations, the tiling conditions on the digit sets play an important role. An interesting question arises naturally:

 \medskip

 \noindent{\bf Question:} {\it What kind of measures are spectral measures and how are they related to translational tilings?}

 \medskip

This question seems to be out of reach using our current knowledge. In this paper,
 we aim to describe a unifying framework bridging the gap between singular spectral
 measures and spectral sets. Let us introduce some simple notations.
 Denote by ${\mathcal L}$ the Lebesgue measure in ${\Bbb R}^d$ and by
 ${\mathcal L}_E$  the normalized Lebesgue measure restricted to the measurable set $E$
(i.e.~${\mathcal L}_E(F) = {\mathcal L}(E\cap F)/{\mathcal L}(E)$). For a finite set $A$,  we
denote by $|A|$ the cardinality of $A$ and by $\delta_A$ the measure $\sum_{a\in A}\delta_a$,
where $\delta_a$ is the Dirac mass at $a$.  We also write
$A\oplus B = C$ if every element in $C$  can be uniquely expressed as
a sum $a+b$ with $a\in A$ and $b\in B$. We now make some observations about
specific examples of  spectral measures known in the literature.

\bigskip

(1)  According to \cite{[PW]}, if $A\subset{\Bbb Z}$ and the set
$\Omega = A+[0,1)$  tiles $[0,N)$, then $\Omega$ is a spectral set.
 We can thus find a set $B$ such that $A\oplus B = \{0,1,...,N-1\}$.
 This means that $\left(\frac{1}{|B|}\delta_{B}\right)\ast {\mathcal L}_{\Omega} = {\mathcal L}_{[0,N]}$.

\medskip

(2) Let $\mu$ be the standard one-fourth Cantor (probability)
measure defined by the self-similar identity
 $$
 \mu(\cdot) = \frac{1}{2}\mu(4\cdot)+ \frac{1}{2}\mu(4\cdot-2).
 $$
 It is known that $\mu$ is a spectral measure \cite{[JP]}. At the same time,
we observe that if we define  $\nu$ to be the one-fourth Cantor measure obeying the equation
$$
\nu(\cdot) = \frac{1}{2}\nu(4\cdot)+ \frac{1}{2}\nu(4\cdot-1),
$$
then $\mu\ast\nu = {\mathcal L}_{[0,1]}$. This can be seen directly by
computing the Fourier transform of both measures.

\bigskip

In fact, we may view the operation of convolution with a positive measure as certain
kind of generalized translation. The above examples suggest the following question.
Let $Q$ be a fundamental domain of some full-rank lattice on ${\Bbb R}^d$.

\medskip

{\bf ${\bf {\mathcal F}(Q)}$:} {\it Any positive Borel measures $\mu$ and $\nu$
such that $\mu\ast\nu = {\mathcal L}_Q$ are spectral measures.}

\medskip

Unfortunately, we cannot expect the above statement
to be  true for all $Q$. In fact, if $\mu = {\mathcal L}_{E}$ with
$E$ is the translational tile without a spectrum constructed in
\cite{[KM]}, then $\mu\ast\nu = {\mathcal L}_Q$ for some
fundamental domain $Q$ as seen directly from the construction of
this counterexample. However, in order to understand which measures are spectral,
it is useful to know  to what extent the statement ${\bf {\mathcal F}(Q)}$
is true for some specific $Q$.  Our first main result unifies the examples
of discrete spectral measures, spectral sets and the singular spectral measures given in (1) and (2) above.

\begin{theorem}\label{th0.1}
For any $d\geq 1$, the statement ${\mathcal F}([0,1]^d)$ is true. Moreover, for any positive Borel
measures $\mu$ and $\nu$ such that  $\mu\ast\nu = {\mathcal L}_{[0,1]^d}$,
we can find spectra $\Lambda_{\mu}$ and $\Lambda_{\nu}$ for $\mu$ and $\nu$
respectively satisfying the property that
$$
\Lambda_{\mu}\oplus\Lambda_{\nu} = {\Bbb Z}.
$$
\end{theorem}

\medskip

We now give a brief explanation of the proof of Theorem \ref{th0.1}. We first focus on
${\Bbb R}^1$ where the proof involves two main steps. The first step is a complete characterization of the
  Borel probability measures $\mu$ and $\nu$ satisfying the identity $\mu\ast\nu = {\mathcal L}_{[0,1]}$.
This characterization is actually a known result in probability due to Lewis \cite{[Le]}. 
 In particular, Lewis proved that only two cases could occur: either one measure is
absolutely continuous and the other one is purely discrete or they are both singular.
To prove our theorem, we will express the measures $\mu$ and $\nu$ as weak limits of
convolutions of some discrete measures using the result of Lewis (See Section 2).
The second step is to construct spectra for $\mu$ and $\nu$. This is done by
observing that the discrete measures obtained at each level are
spectral measures. We then show that the spectral property carries over by
passing to the weak limit. This argument is a generalization of the proof in
\cite{[DHL]} (See Section 3). After the dimension one case is established,
we characterize the Borel probability measures $\mu$ and $\nu$ satisfying
$\mu\ast\nu = {\mathcal L}_{[0,1]^d}$ as Cartesian products of one-dimensional
Borel probability measures $\sigma_i$ and $\tau_i$,
$i=1,...,d$, on ${\Bbb R}^1$
satisfying $\sigma_i\ast\tau_i={\mathcal L}_{[0,1]}$
and also prove  the spectral property for those (See Section 5).

\medskip

It is very unclear whether ${\bf {\mathcal F}({Q})}$ is true if $Q$ is not a hypercube.
We will focus our attention on ${\Bbb R}^1$ in which Fuglede's conjecture
remains open. We propose the
following generalized Fuglede's conjecture for spectral measures on
${\Bbb R}^1$ and it is direct to see that a full generality of
${\bf {\mathcal F}({Q})}$ on ${\Bbb R}^1$ will imply one direction of this generalized conjecture.

\medskip

\noindent{\bf Conjecture (Generalized Fuglede's Conjecture):}
{\it A compactly supported Borel probability measure $\mu$ on ${\Bbb R}^1$ is
spectral if and only if there exists a Borel probability measure
$\nu$ and a fundamental domain $Q$ of some lattice  on ${\Bbb R}^1$
such that $\mu\ast\nu ={\mathcal L}_{Q}$.}

\medskip

This is an open conjecture on ${\Bbb R}^1$ and we will prove that
it extends the classical Fuglede's conjecture.

\medskip

\begin{theorem}
The generalized Fuglede's conjecture implies Fuglede's conjecture on ${\Bbb R}^1$.
\end{theorem}

 Let us make some remarks on the classical Fuglede's conjecture on ${\Bbb R}^1$.
There is some evidence that the conjecture may be true on ${\Bbb R}^1$.
 In particular, the known fact that all tiling sets of a tile and all spectra of a
spectral set are  periodic offers some credibility to the conjecture \cite{[LW1],[IK]}.
Moreover, some algebraic conditions, if satisfied, are sufficient to settle
the conjecture on ${\Bbb R}$, although these conditions are not easy to check \cite{[DL2]}.

%
%

\medskip

As our focus is the one-dimensional case, we organize our paper as follows:
In Section 2, we describe the factorization of the Lebesgue measure on
$[0,1]$ given by Lewis and, for the reader's convenience, we provide a somewhat
different proof of the factorization theorem that avoids some of the complications
of the original ones stemming from the use of probabilistic tools.
We then prove the spectral property in Section 3 and discuss the
generalized Fuglede's conjecture on ${\Bbb R}^1$ in Section 4.
We will finally prove Theorem \ref{th0.1} in higher dimension in Section 5.
As this piece of work offers us several new directions for further research, we end this paper with some remarks and open question in Section 6.

\medskip

\noindent{\it Note:} During the preparation of the manuscript,
we were made aware that Professor Xinggang He  and his student \cite{[AH]}
discovered independently a new class of one-dimensional spectral measures
obtained via a Moran construction of fractals. These
one-dimensional spectral measures turn out to coincide exactly
with those we consider in this paper.

\bigskip

\section{Factorization of Lebesgue measures}

%

Let ${\mathcal L}_{[0,1]}$ be the Lebesgue measure supported on $[0,1]$ and
let $\mu$ and $\nu$ be two Borel probability measures supported on $[0,1]$.
We say that $(\mu,\nu)$ is a {\it complementary pair} of measures with respect to
${\mathcal L}_{[0,1]}$ if
$$
\mu\ast\nu ={\mathcal L}_{[0,1]}.
$$
Let ${\mathcal N} = \{N_k\}_{k=1}^{\infty}$ be a sequence of positive
integers greater than or equal to 2. We associate  with ${\mathcal N}$ the discrete measures
\begin{equation}\label{eqdis}
\nu_{k}  = \frac{1}{N_k}\sum_{j=0}^{N_k-1}\delta_{\frac{j}{N_1\cdots N_k}},\quad k\ge 1.
\end{equation}
 For a given Borel set $E$, recall that ${\mathcal L}_{E}$ is the normalized
Lebesgue measure supported on $E$. We now  observe that the Lebesgue measure supported on
$[0,1]$ admits a natural decomposition as convolution products.
$$
\begin{aligned}
{\mathcal L}
_{[0,1)} =& \nu_{1}\ast ({\mathcal L}_{[0,\frac{1}{N_1}]})\\
=&\nu_{1}\ast\nu_{2}\ast({\mathcal L}_{[0,\frac{1}{N_1N_2}]})\\
=&\cdots\\
 =& \nu_{1}\ast\nu_{2}\ast\cdots\ast\nu_{k}\ast({\mathcal L}_{[0,\frac{1}{N_1\cdots N_k}]}).
\end{aligned}
$$
The sequence of measures $\nu_{1}\ast\nu_{2}\ast\cdots\ast\nu_{k}$ converges weakly to
${\mathcal L}_{[0,1]} $. Therefore, one can write the Lebesgue measure as an infinite
convolution of discrete measures.
\begin{equation}\label{eq0.1}
{\mathcal L}_{[0,1]} = \nu_{1}\ast\nu_{2}\ast\cdots.
\end{equation}

\medskip

 Given a set ${\mathcal N}$ as above, we will consider two types of factorization (Type I and Type II)
 of ${\mathcal L}_{[0,1]}$  as the convolution of two measures obtained from the
infinite factorization obtained in
(\ref{eq0.1}).

\medskip

   \noindent{\bf Type I.} There exists a finite positive integer  $k$ such that
we have either
   $$
   \mu_{\mathcal N} = \nu_{1}\ast\nu_{3}\ast...\ast\nu_{{2k-1}} \  \mbox{and} \
   \nu_{\mathcal N} = \nu_{2}\ast\nu_{4}\ast...\ast\nu_{{2k}}\ast({\mathcal L}_{[0,\frac{1}{N_1N_2\cdots
N_{2k}}]})
   $$
 or
   $$
   \mu_{\mathcal N} = \nu_{1}\ast\nu_{3}\ast...\ast\nu_{{2k-1}}\ast({\mathcal L}_{[0,\frac{1}{N_1N_2
\cdots N_{2k}}]}) \  \mbox{and} \ \nu_{\mathcal N} = \nu_{2}\ast\nu_{4}\ast...\ast\nu_{{2k}}.
   $$

\medskip

\noindent{\bf Type II}
\begin{equation}\label{eq0.2}
 \mu_{{\mathcal N}} = \nu_{1}\ast\nu_{3}\ast\cdots\ast\nu_{{2k-1}}\ast\cdots
\end{equation}
\begin{equation}\label{eq0.3}
\nu_{{\mathcal N}} = \nu_{2}\ast\nu_{4}\ast\cdots\ast\nu_{{2k}}\cdots
\end{equation}

\medskip

\begin{Rem}\label{rem1}
 {\rm The reader might want to construct more general decompositions obtained by choosing other
factorizations of (\ref{eq0.1}), but note that if convolution product of two consecutive factors of
 (\ref{eq0.1}) belong to the same factor in the factorization, say $\nu_k$ and $\nu_{k+1}$, then we have}
$$
\nu_k\ast\nu_{k+1} = \frac{1}{N_kN_{k+1}}\sum_{j=0}^{N_kN_{k+1}}\delta_{j/N_1N_2...(N_kN_{k+1})}
$$
{\rm and we would then be able to write the given convolution product as one of type I or type II associated with a different ${\mathcal N}$.}
\end{Rem}
\medskip

Note in both cases that
$\mu_{{\mathcal N}}\ast\nu_{{\mathcal N}} ={\mathcal L}_{[0,1]} $ by (\ref{eq0.1}).
Therefore, they are $\mu_{{\mathcal N}}$ and $\nu_{{\mathcal N}}$
form a complementary pair with respect to ${\mathcal L}_{[0,1]}$.
In the case of the Type I decomposition, one is  purely discrete and one is absolutely continuous
while in the Type II decomposition,  both factors are singularly continuous measures.
We say that a complementary pair $(\mu, \nu)$ is {\it natural} if
we can find a sequence ${\mathcal N}$ of positive integers such that
$(\mu, \nu) =(\mu_{{\mathcal N}},\nu_{{\mathcal N}})$.

\begin{theorem}\label{th1.1}
If $\mu$ and $\nu$ are positive Borel probability measures supported on $[0,1]$
 and $\mu\ast \nu = {\mathcal L}_{[0,1]}$, then $\mu$ and $\nu$ are natural complementary pair.
\end{theorem}

\medskip

This theorem is essentially due to Lewis \cite{[Le]} who considered the problem
in probability consisting in characterizing the type of the distributions of pairs of
independent random variables $X$ and $Y$ whose sum $X+Y$ is a uniform random variable on
$[-\pi,\pi]$. For the reader's convenience, we will give here another proof based on his ideas
as his result is not widely known. Moreover, the proof we give here is more analytical in flavor and
 avoids some of the complications arising in the original proof from the use of probability tools.
 The main important step of the proof is to show that if
two probablity measures $\mu$ and $\nu$ satisfy $\mu\ast\nu = {\mathcal L}_{[0,1]}$,
then one of them, say $\mu$, must be ''$1/N$ periodic" in the sense that
$\mu = \left(1/N\sum_{j=0}^{N-1}\delta_{j/N}\right)\ast \mu_1$ for some integer
$N\geq 2$ and $\mu_1\ast\nu = {\mathcal L}_{[0,1/N]}$. This is done by analyzing the
structure of the zeros of the Fourier transform of $\mu$ and $\nu$ (Lemma \ref{lem1}).

\medskip

We now define the (complex) Fourier transform of a compactly supported probability measure $\mu$ by the formula
$$
\widehat{\mu}(\xi) = \int e^{-2\pi i \xi x}d\mu(x), \ \xi\in{\Bbb C}.
$$
We will consider convolution products yielding the Lebesgue measure supported on $[-1/2,1/2]$ instead of $[0,1]$ to exploit some symmetric properties of the solutions (as explained below).  Note that $\mu\ast \nu = {\mathcal L}_{[-1/2,1/2]}$ is equivalent to
\begin{equation}\label{sin}
\widehat{\mu}(\xi)\widehat{\nu}(\xi)=\widehat{{\mathcal L}_{[-1/2,1/2]}}(\xi) =\frac{\sin \pi \xi}{\pi \xi}.
\end{equation}
The zero set of the Fourier transform $\widehat{\mu}$ in the complex plane will be denoted by
 $$
 {\mathcal Z}(\widehat{\mu}) = \{\xi\in {\mathbb C}: \widehat{\mu}(\xi) =0\}
 $$
 Since $((\delta_x\ast\mu)\ast (\delta_{-x}\ast\nu) = {\mathcal L}_{[-1/2,1/2]}$
for any real numbers $x$, we may assume the smallest closed interval containing the support of
$\mu$ is given by $[-a,a]$. Denote by supp $\mu$ the closed support of $\mu$. Given a probability measure
$\rho$, we also define the measure $\check{\rho}$ to be the measure satisfying
$\check{\rho}(B) = \rho(-B)$ for any Borel set $B\subset \mathbb{R}$.

 \begin{Lem}\label{lem01}
Let $\mu$ and $\nu$ be two probability measures such that $\mu\ast\nu = {\mathcal L}_{[-1/2,1/2]}$ and assume that the smallest closed interval containing supp $\mu$ is of the form $[-a,a]$, $a>0$. Then we have
\begin{equation}\label{eq1.1}
{\mathbb Z}\setminus\{0\} = {\mathcal Z}(\widehat{\mu})\cup {\mathcal Z}(\widehat{\nu}) \ (\mbox{\rm as a disjoint union}).
\end{equation}
Moreover, the smallest closed interval containing  supp $\nu$ is given by $[-b,b]$ where $b=1/2-a$ and both  $\mu$ and $\nu$ have symmetric distributions around the origin (i.e. $\check{\mu} = \mu$ and $\check{\nu} = \nu$).
 \end{Lem}

 \begin{pf}
It is well-known that $\widehat{\mu}$ is a non-zero entire analytic function, so its zero set is
a discrete set in the complex plane. Furthermore, since the zeros of $\widehat{\chi_{[-1/2,1/2]}}$ are
simple, (\ref{eq1.1}) follows from (\ref{sin}). Let $[c,b]$ be the smallest closed interval
containing the support of $\nu$. Then $a+b=1/2$ and $-a+c=-1/2$ showing that $c=-b$ and $b=1/2-a$.
\medskip

Finally, note that, since $\mu$ is a positive measure,
${\mathcal Z}\left(({\check{\mu}})^{\widehat{}}\right) =
{\mathcal Z}({\widehat{\mu}})$. Therefore,
${\mathcal Z}\left(({\check{\mu}})^{\widehat{}}\right)\cup{\mathcal Z}(\widehat{\nu})
= {\Bbb Z}\setminus\{0\}$. Consider the tempered distribution
$\rho: = \check{\mu}\ast\nu\ast\delta_{\Bbb Z}$.
Then $\widehat{\rho} =({\check{\mu}})^{\widehat{}}\cdot\widehat{\nu}\cdot\delta_{\Bbb Z}
= \delta_0$. Hence, $\rho$ is the Lebesgue measure on ${\Bbb R}$ and the restriction of $\rho$
to the interval $[-1/2,1/2]$ is $\check{\mu}\ast\nu$. This shows that
$\check{\mu}\ast\nu = {\mathcal L}_{[-1/2,1/2]}$, which means that
$\check{\mu}\ast\nu = \mu\ast\nu$. Taking Fourier transform, we obtain $\check{\mu} = \mu$.
The proof of the symmetry of $\nu$ is similar.
\end{pf}

Note that Lewis used the Hadamard factorization theorem to prove the
symmetry property of $\mu$ and $\nu$ in Lemma \ref{lem01}. The ideas of the following
two lemmas are due to Lewis and form the crucial parts of the argument.

\medskip

\begin{Lem}\label{lem0}
Let $r\geq1$ be the smallest positive zero of $\widehat{\mu}$. Then
$$
\frac{1}{4r}\leq a\leq \frac{1}{2r} \ \mbox{and} \ \frac{1}{2}-\frac{1}{2r}\leq b\leq \frac{1}{2}-\frac{1}{4r}.
$$
\end{Lem}

\medskip

\begin{pf}
We just need to prove the lower estimates for both $a$ and $b$ as the upper ones will follow from these and the fact that $a+b=1/2$. Since $r$ is a zero of $\widehat{\mu}$, then $-r$ is also a zero and we must have $\int \cos (2\pi rx)d\mu(x) =0$. This implies that $2\pi r a \geq \frac{\pi}{2}$ and thus $a\geq\frac{1}{4r}$. In particular, the claim is true for $r=1$.

\medskip

For the upper bound, we consider the following functions for different $r$.
$$
h(x) : = \left\{
           \begin{array}{ll}
             \cos(2\pi x), & \hbox{$r=2$;} \\
             \cos(2\pi x)-\cos(2\pi 2x), & \hbox{$r=3$;} \\
             \cos(\frac{\pi rx}{2})\prod_{j=1}^{k-1}(\cos(2\pi x)-\cos\frac{2(2j-1)\pi}{r}), & \hbox{$r>2$, $r =2k$;} \\
             (\cos(\frac{\pi(r-1)x}{2})-\cos(\frac{\pi(r+1)x}{2})\prod_{j=1}^{k-2}(\cos(2\pi x)-\cos\frac{(2\pi)(2j)}{r}), & \hbox{$r>2$, $r=2k-1$.}
           \end{array}
         \right.
$$
By expanding $h(x)$, we see that $h(x)$ is a linear combination of
 $\cos(2\pi kx)$, for $k =1...,r-1$. Hence $\int h(x)d\nu(x) =0$ as $1,\cdots,r-1$
are zeros of $\widehat{\nu}$. By checking the sign of each factor, we see that
if $2\pi x\leq \pi(r-1)/r$, then  $h(x)\geq 0$.

\medskip

 Consider the case where $r>2$ is even.
 We have either $2\pi b\geq \pi(r-1)/r$ (i.e. $b\geq1/2-1/2r$) or $\nu$ is supported on
 the atoms $\pm (1/r),\cdots, \pm (r-3)/r$. However, $\nu$ cannot be supported on those atoms
since $\widehat{\nu}$ would be a polynomial in $\cos(2\pi x/r)$ of degree at most $r-3$,
but there are $r-1$ zeros for $\widehat{\nu}$, a contradiction. Therefore, we must have $b\geq 1/2-1/(2r)$.
 The proof for the other cases follows from a similar argument.
\end{pf}

\medskip

\begin{Lem}\label{lem1}
Let $N>0$ be a positive integer and let $\mu$ and $\nu$ be two probability measures
on ${\Bbb R}$ such that  $\mu\ast\nu ={\mathcal L}_{[0,1/N]}$ with
 neither $\widehat{\mu}$ nor $\widehat{\nu}$ being identically one.
Suppose that $N\in {\mathcal Z}(\widehat{\nu})$ and let $Nr$ with $r>1 $
be the smallest positive zero of $\widehat{\mu}$. Then
$${\mathcal Z}(\widehat{\mu})\subset Nr{\Bbb Z}.$$ 
\end{Lem}
\begin{pf}
By  rescaling the measures by a factor of $N$, it is easy to see that it suffices
to consider the case $N=1$. By translating the measure
(i.e. $\mu\ast(\delta_{-1/2}\ast\nu) = {\mathcal L}_{[-1/2,1/2]}$),
it suffices to prove the lemma for the case $\mu\ast\nu ={\mathcal L}_{[-1/2,1/2]}$,
where $\check{\mu} = \mu$ and $\check{\nu} = \nu$.

\medskip

Let $\rho (E) = \nu(\{0\})\delta_0(E)+2\nu(E\cap (0,1/2])$ and $\check{\rho}(E) = \rho(-E)$ for $E$ Borel. Then,  the fact that $\nu(E) = \nu(-E)$ implies that $
\rho+\check{\rho}= 2\nu.$
Therefore,
\begin{equation}\label{eqpf1.2}
\mu\ast\rho+\mu\ast\check{\rho} =2{\mathcal L}_{[-1/2,1/2]}.
\end{equation}
This implies, in particular, that $\mu\ast\rho$ is absolutely continuous
with respect to the Lebesgue measure
and we can let $g(x)\geq0$ be its density. Then $g(-x)$ is the density of
$(\mu\ast\rho)^{\check{}} = \mu\ast\check{\rho}$. By (\ref{eqpf1.2}),
$$
g(x)+g(-x) =2, \ a.e.
$$
As supp ($\mu\ast\check{\rho}$) (and hence supp $g(-x)$) is contained in $[-1/2,a]$, $g(x)=2$  on $[a,1/2]$. We may therefore write
$$
\begin{aligned}
g = 2\chi_{[a,1/2]}+g\chi_{[-a,a]} =& 2\chi_{[a,1/2]}+ g\chi_{[-a,0]} +(2-g(-x))\chi_{[0,a]} \\
=& 2\chi_{[0,1/2]}+( g\chi_{[-a,0]} - g(-x)\chi_{[0,a]}).
\end{aligned}
$$
Note that $2\chi_{[0,1/2]}$ is the density of the measure ${\mathcal L}_{[0,1/2]}$. Taking Fourier transform, we have
\begin{equation}\label{eqpf1.3}
\widehat{\mu}(\xi)\widehat{\rho}(\xi) = \widehat{g}(\xi) = \widehat{{\mathcal L}_{[0,1/2]}}(\xi)+2i\int_{0}^{a} g(-x)\sin(2\pi \xi x)dx
\end{equation}
Suppose that $r$ is even. As $\widehat{\mu}(r)=0$, we must have
$$
\int_{0}^{a} g(-x)\sin(2\pi r x)dx=0.
$$
Since $a\leq 1/2r$ by Lemma \ref{lem0}, we have $\sin(2\pi rx)\geq 0$  on $[0,a]$ and thus $g(-x)=0$ there. Thus, (\ref{eqpf1.3}) implies that
\begin{equation}\label{eqpf1.4}
\widehat{\mu}(\xi)\widehat{\rho}(\xi)  = \widehat{{\mathcal L}_{[0,1/2]}}(\xi).
\end{equation}
Hence, ${\mathcal Z}(\widehat{\mu}) \subset2{\Bbb Z}$.

\medskip

 Writing  $r= 2^nm$ where $m$ is odd, we deduce from the above argument
that ${\mathcal Z}(\widehat{\mu}) \subset2{\Bbb Z}$. Consider the measure $\mu_1(E) = \mu(E/2)$ and $\rho_1(E) = \rho(E/2)$ we have $\widehat{\mu_1}(\xi) = \widehat{\mu}(2\xi)$ and $\widehat{\rho_1}(\xi) = \widehat{\rho}(2\xi)$. By (\ref{eqpf1.4}), we have
$\widehat{\mu_1}(\xi) \widehat{\rho_1}(\xi) =\widehat{{\mathcal L}_{[0,1]}}(\xi)$
(i.e. $\mu_1\ast(\delta_{-1/2}\ast\rho_1) = {\mathcal L}_{[-1/2,1/2]}$).
Moreover, ${\mathcal Z}(\widehat{\mu_1}) = \frac12{\mathcal Z}(\widehat{\mu})$.
In this case, the smallest positive zero of $\widehat{\mu_1}$ will be $2^{n-1}m$. Therefore,
repeating the above argument, we have ${\mathcal Z}(\widehat{\mu})\subset 2^n{\Bbb Z}$
and the proof will be finished if we can prove our claim if $r$ is odd.

\medskip

Suppose now that $r$ is odd.  We consider the measures
$\nu_1(E) = \nu(E\cap [-a,b])$ and $\nu_2(E) = \nu(E\cap [-b,-a))$
(Here, it is more convenient not to normalize $\nu_1$ and $\nu_2$ as probability
measures). We have then  $\nu = \nu_1+\nu_2$ and $
{\mathcal L}_{[-1/2,1/2]} = \mu\ast\nu_1+\mu\ast\nu_2.$
Let $g_1$ and $g_2$ be the density of $\mu\ast\nu_1$ and $\mu\ast\nu_2$ respectively. The above implies that
$$
g_1(x) +g_2(x) = 1 \ \mbox{a.e. on} \ [-1/2,1/2].
$$
Note that the supp $g_1$ is contained in $[-2a,1/2]$ and supp $g_2$ is contained in $[-1/2,0]$. It follows that $g_1 = 1$ almost everywhere on $[0,1/2]$. We may therefore write
$$
g_1 = \chi_{[0,1/2]}+g_1\chi_{[-2a,0]}.
$$
Taking Fourier transforms  and noting that $\widehat{g_1}(\xi)= \widehat{\mu}(\xi)\widehat{\nu_1}(\xi)$, we obtain
\begin{equation}\label{eqpf1.5}
\widehat{\mu}(\xi)\widehat{\nu_1}(\xi) = \widehat{\chi_{[0,1/2]}}(\xi) + \int_{0}^{2a}g_1(-x)e^{2\pi i \xi x}dx.
\end{equation}
As $\widehat{\mu}(r)=0$, by substituting $\xi = r$ and equating the imaginary parts, we have
$$
    \frac{1}{\pi r} = \int_0^{2a}g_1(-x)\sin(2\pi  r x)dx.
$$
By Lemma \ref{lem0}, $2a\geq1/2r$ and therefore,
$$
\begin{aligned}
\frac{1}{\pi r} =& \int_{0}^{1/2r}g_1(-x)\sin(2\pi  r x)dx+\int_{1/2r}^{2a}g_1(-x)\sin(2\pi  r x)dx\\
\leq&\int_{0}^{1/2r}g_1(-x)\sin(2\pi  r x)dx \ \  \ ( \mbox{as} \  \sin(2\pi rx)\leq0 \  \mbox{on} \ [1/2r,2a]) \\
\leq&\int_{0}^{1/2r}\sin(2\pi  r x)dx =\frac{1}{\pi r}.  \  \  \ ( \mbox{as} \ g_1(-x)\leq 1)\\
\end{aligned}
$$
Hence, we must have $g_1(-x)=1$ on $[0,1/2r]$ and $\int_{1/2r}^{2a}g_1(-x)\sin(2\pi  r x)dx=0$, which implies that $g_1(-x)=0$ on $[1/2r,2a]$. Considering the real part of the equation (\ref{eqpf1.5}) and noting that $\widehat{\mu}(\xi)$ is real-valued (as $\check{\mu} = \mu$), we have
$$
\widehat{\mu}(\xi){\mbox Re}\left(\widehat{\nu_2}(\xi)\right) = \frac{\sin\pi\xi}{2\pi\xi}+\int_{0}^{1/2r}\cos(2\pi \xi x)dx = \frac{1}{2\pi \xi}\left(\sin\pi\xi+\sin\frac{\pi\xi}{r}\right).
$$
Since ${\mathcal Z}(\widehat{\mu})\subset {\Bbb Z}$, the previous equation shows that in fact  ${\mathcal Z}(\widehat{\mu})\subset r{\Bbb Z}$, completing the proof.
\end{pf}

\medskip

\noindent{\it Proof of Theorem \ref{th1.1}.}
Let $(\mu,\nu)$ be a complementary pair with respect to ${\mathcal L }_{[0,1]}$.
We may assume that $\widehat{\nu}(1)\neq 0$ and we let $N_1>1$ be
the smallest positive zero of $\widehat{\nu}$. We have ${\mathcal Z}(\widehat{\nu})\subset N_1{\Bbb Z}$
by Lemma \ref{lem1}. As the zero sets of $\widehat{\mu}$ and $\widehat{\nu}$ are disjoint
 (see (\ref{eq1.1})), the set $\{k\in{\Bbb Z}: \widehat{\mu}(k)\neq 0\}$ is contained in $N_1{\Bbb Z}$.

\medskip

Consider the periodization of the measure $\mu$ defined by $\mu_p = \mu\ast\delta_{\Bbb Z}$.
Its distributional Fourier transform (as a tempered distribution) is given by
$$
\widehat{\mu_p} = \widehat{\mu}\cdot\delta_{\Bbb Z} =\widehat{\mu}\cdot\delta_{N_1{\Bbb Z}}
$$
Hence, $\mu_p$ is indeed $1/N_1$-periodic. It follows immediately that
\begin{equation}\label{eq1.2}
\mu =\nu_1 \ast\alpha_1 \ \  \ \mbox{and}  \ \ \ \nu\ast\alpha_1 ={\mathcal L}_{[0,1/N_1]}
\end{equation}
where $\nu_1  =\frac{1}{N_1}\sum_{j=0}^{N_1-1}\delta_{j/N_1}$ and
$\alpha_1(E) = N_1\mu(E\cap{[0,1/N_1]})$ for any Borel set $E$. The case where $\alpha_1$ is the
Dirac measure at the origin  immediately yields a type I decomposition.
Otherwise, we apply Lemma \ref{lem1} on the pair $(\nu,\alpha_1)$.
Since $\widehat{\nu}(N_1)=0$, we have $\widehat{\alpha_1}(N_1)\neq0$ and we
can let $N_2$ ne the smallest positive integer such that $\widehat{\alpha_1}(N_1N_2)=0$.
By Lemma \ref{lem1}, we have ${\mathcal Z}(\widehat{\alpha_1})\subset N_1N_2{\Bbb Z}$.
We obtain
$$
\mu = \nu_1\ast\alpha_1, \ \ \ \nu = \nu_2\ast\alpha_2 \ \ \ \alpha_1\ast\alpha_2 = {\mathcal L}_{[0,1/N_1N_2]}
$$
where $\nu_2 = \frac{1}{N_2}\sum_{j=0}^{N_2-1}\delta_{j/N_1N_2}$.
The case where $\alpha_2$ is a Dirac measure at the origin yelds again a type I decomposition.
Otherwise, we continue this inductive process and define recursively the probability measures
$\alpha_k$, $k\geq 1$. If $\alpha_k=\delta_0$ for some $k$, the process stops
and we have arrived at a type I decomposition. If $\alpha_k\neq\delta_0$ for all $k$, we have
then expressed both measures $\mu$ and $\nu$ at the infinite convolution products
$\mu = \nu_1\ast\nu_3\ast \dots,  \quad \nu = \nu_2\ast\nu_4\ast \dots$,
which yields a type II decomposition.
\qquad$\Box$

\medskip

Theorem \ref{th1.1} also gives us a new proof of classification
of the set $A$ and $B$ such that  $A\oplus B = \{0,...,n-1\}$ which was proved in
\cite{[Lo]} and \cite{[PW]} using a theorem of De Bruijn.

\medskip

\begin{Cor}\label{cor2.1}
Let  ${\mathcal E}_n = \{0,1,\cdots,n-1\}$ and let
 ${\mathcal A}$ and ${\mathcal B}$ be two finite set of integers such that
${\mathcal A}\oplus {\mathcal B} = \{0,...,n-1\}$. Suppose that $1\in {\mathcal A}$.
Then there exist integers $N_1,...,N_{2k}$ such that $N_1...N_{2k} =n$ and
$$
A = {\mathcal E}_{N_0}\oplus N_0N_1{\mathcal E}_{N_2}\oplus...\oplus
 N_0N_1...N_{2k-1}{\mathcal E}_{2k}
$$
$$
B = N_0{\mathcal E}_{N_1}\oplus N_0N_1N_2{\mathcal E}_{N_3}\oplus...\oplus
N_0N_1...N_{2k-2}{\mathcal E}_{2k-1}.
$$
\end{Cor}

\medskip

\begin{pf}
As   ${\mathcal A}\oplus {\mathcal B} = \{0,...,n-1\}$, we have
$$
\left(\frac{1}{|{\mathcal A}|}\delta_{\frac{1}{n}{\mathcal A}}\right)\ast
\left(\frac{1}{|{\mathcal B}|}\delta_{\frac{1}{n}{\mathcal B}}\right)\ast
{\mathcal L}_{[0,1/n]} = {\mathcal L}_{[0,1]}.
$$
By Theorem \ref{th1.1}, the measures
$\mu=\left(\frac{1}{|{\mathcal A}|}\delta_{\frac{1}{n}{\mathcal A}}\right)$
and $\nu =\left(\frac{1}{|{\mathcal B}|}\delta_{\frac{1}{n}{\mathcal B}}\right)\ast
{\mathcal L}_{[0,1/n]}$ are natural complementary pair.
As one of them is discrete and the other is absolutely continuous, they correspond to a type I
decomposition. Since $1\in{\mathcal A}$, we have thus $1/n\in \frac{1}{n}{\mathcal A}$.
By comparing the support of the measures, we obtain the existence of integers $N_1',N_2'...$ such that
$$
\frac1n{\mathcal A} = \frac{1}{N_1'}{\mathcal E}_{N_1'}\oplus
\frac{1}{N_1'N_2'N_3'}{\mathcal E}_{N_3'}\oplus...\oplus
\frac{1}{N_1'N_2'...N_{2k-1}'}{\mathcal E}_{N_{2k-1}'}.
$$
$$
\frac1n{\mathcal B} = \frac{1}{N_1'N_2'}{\mathcal E}_{N_2'}\oplus
\frac{1}{N_1'N_2'N_3'N_4'}{\mathcal E}_{N_4'}\oplus...\oplus
\frac{1}{N_1'N_2'...N_{2k}'}{\mathcal E}_{N_{2k}'}
$$
and $n = N_1'...N_{2k}'$. Letting $N_r = N'_{2k-r}$, we obtain the desired factorization.
\end{pf}

\bigskip

\section{The spectral property}

In this section, we show that all measures appearing in natural complementary pairs
are spectral measures. Recall that a Borel probability measure $\mu$ is called a {\it spectral measure} with associated {\it spectrum} $\Lambda$ if the collection of exponentials $E(\Lambda) = \{e^{2\pi i \lambda x}\}_{\lambda\in\Lambda}$ forms an orthonormal basis for $L^2(\mu)$. It is easy to see that $E(\Lambda)$ is an orthonormal set in $L^2(\mu)$ if and only if
  $$
  \Lambda-\Lambda\subset {\mathcal Z}(\widehat{\mu})\cup\{0\}.
  $$
  By a well-known result in \cite{[JP]}, $\Lambda$ is a spectrum of $\mu$ if and only if
  \begin{equation}\label{eq4.0}
  Q(\xi): = \sum_{\lambda\in\Lambda}|\widehat{\mu}(\xi+\lambda)|^2\equiv 1.
  \end{equation}
  In fact, if $E(\Lambda)$ is an orthonormal set,
$Q(\xi)\leq 1$ and $Q$ is an entire function of exponential type (\cite{[JP]}, see also \cite{[DHL]}).
Let ${\mathcal N} = \{N_i\}_{i=1}^{\infty}$ be a collection of positive integers and
consider the Type I and II decomposition as in the previous section. Let
 $$
 \mu^{(k)} = \nu_{1}\ast\nu_{3}\ast\cdots\ast\nu_{2k-1}, \
\nu^{(k)} = \nu_{2}\ast\nu_{4}\ast\cdots\ast\nu_{2k}
 $$
 and for a given ${\mathcal N}$, we let $A_1 = \{0,..,N_1-1\}$ and
$A_n = N_1\cdots N_{n-1}\cdot\{0,..,N_n-1\}$ for $n\geq2$. We start with a simple observation.

 \medskip
 \begin{Prop}\label{prop4.1}
Each $\nu_{n}$ is a spectral measure with spectrum $A_n$.
 For all $k\geq 1$, $\mu^{(k)}$ is a spectral measure with spectrum given by
\begin{equation}\label{eq4.1}
 \Lambda_k = \bigoplus_{j=1}^{k}A_{2j-1}
\end{equation}
 In particular, the type I natural complementary pair $\mu_{\mathcal N}$ and $\nu_{\mathcal N}$ defined in the previous section are spectral measures.
 \end{Prop}

 \begin{pf}
It is immediate to see that the measure $\frac{1}{N_n}\sum_{j=0}^{N_n-1}\delta_{j/N_n} $ is a
 spectral measure with spectrum $\{0,..,N_n-1\}$.
Therefore, $\nu_{n} = \frac{1}{N_n}\sum_{j=0}^{N_n-1}\delta_{j/(N_1\cdots N_{n})}$ is a
spectral measure with spectrum  $N_1\cdots N_{n-1}\cdot\{0,..,N_n-1\} =A_n$.

\medskip

Note that ${\mathcal Z}(\widehat{\nu_n}) = N_1 N_2...N_n{\Bbb Z}\setminus
N_1 N_2...N_{n-1}{\Bbb Z}$ and
  $$
\widehat{\mu^{(k)}}(\xi) = \prod_{j=1}^{k}\widehat{\nu_{2j-1}}(\xi).
$$
 For notational convenience, we define $N_0=1$.
Taking distinct $\lambda_1,\lambda_2\in\Lambda_k$ and writing
$\lambda_{\ell} = \sum_{j=1}^{k}r_{\ell,j}N_1N_2...N_{2j-2}$, for $\ell=1,2$, we have
 $$
\lambda_1-\lambda_2 = \sum_{j=1}^{k}(r_{1,j}-r_{2,j})N_1N_2\cdots
 N_{2j-2}=\sum_{j=J}^{k}s_jN_1N_2\cdots N_{2j-2},
$$
where $J$ is the first index such that $r_{1,j}\neq r_{2,j}$ and
$-(N_{2J-1}-1)\leq s_J\leq N_{2J-1}-1$. so
$\widehat{\nu_{2J-1}}(\lambda_1-\lambda_2) = \widehat{\nu_{2J-1}}(N_1\cdots N_{2J-2}s_J)=0$.
Therefore, $\widehat{\mu^{(k)}}(\lambda_1-\lambda_2)=0$.
This proves the orthogonality of $E(\Lambda_k)$ in  $L^2(\mu^{(k)})$.
As $L^2(\mu^{(k)})$ is a finite dimensional vector space of dimension
$N_{1} N_3\cdots N_{2k-1}=\mbox{card}\left(E(\Lambda)\right)$, the collection $E(\Lambda)$
 must be complete  in  $L^2(\mu^{(k)})$.

\medskip

To prove the last statement, we just consider the case where $\mu_{\mathcal N} = \mu^{(k)}$
and  $\nu_{\mathcal N} = \nu^{(k)} \ast ({\mathcal L}_{[0,\frac{1}{N_1\cdots N_{2k}}]})$,
as the case $
   \mu_{\mathcal N} = \nu_{1}\ast\nu_{3}\ast...\ast\nu_{{2k-1}}\ast({\mathcal L}_{[0,\frac{1}{N_1N_2
\cdots N_{2k}}]})$ and $\nu_{\mathcal N} = \nu_{2}\ast\nu_{4}\ast...\ast\nu_{{2k}}$
is similar.
It is easily seen, as before, that $\nu^{(k)}$ is also a discrete spectral measure with spectrum
$$
\widetilde{\Lambda_{k}} = \bigoplus_{j=1}^{k}A_{2j}.
$$
 Moreover, $\widehat{\nu^{(k)}}$ is $N_1...N_{2k}$-periodic.
Let $\alpha$ denote the measure ${\mathcal L}_{[0,\frac{1}{N_1N_2\cdots N_{2k}}]}$. Then $\alpha$ has $N_1N_2...N_{2k}{\Bbb Z}$ as a spectrum. It follows that
$$
\begin{aligned}
\sum_{\lambda\in\widetilde{\Lambda_{k}}+N_1\cdots N_{2k}{\Bbb Z}}|\widehat{\nu_{\mathcal N}}(\xi+\lambda)|^2=&\small{\sum_{\lambda\in\widetilde{\Lambda_{k}},m\in{\Bbb Z}}}|\widehat{\nu^{(k)}}(\xi+\lambda+N_1...N_{2k}m)|^2|\widehat{\alpha}(\xi+\lambda+N_1...N_{2k}m)|^2\\
=&\sum_{\lambda\in\widetilde{\Lambda_{k}}}|\widehat{\nu^{(k)}}(\xi+\lambda)|^2\cdot\sum_{m\in{\Bbb Z}}|\widehat{\alpha}(\xi+\lambda+N_1...N_{2k}m)|^2\equiv1.\\
\end{aligned}
$$
Hence, $\nu_{\mathcal N}$  a spectral measure with
spectrum $\widetilde{\Lambda_{k}}+N_1\cdots N_{2k}{\Bbb Z}$.
 \end{pf}

 \medskip

 It remains to deal with the spectral property for complementary pairs  $\mu_{\mathcal N}$ and
$\nu_{\mathcal N}$ of type II. Since these two measures have essentially the same form,
we will discuss only the case $\mu: = \mu_{\mathcal N}$. Note that
the measure $\mu$ will be the weak limit of the measures $\mu^{(k)}$ and
\begin{equation}\label{eq4.4}
\widehat{\mu}(\xi)= \prod_{j=1}^{\infty}\widehat{\nu_{2j-1}}(\xi) =
\widehat{\mu^{(k)}}(\xi)\cdot\prod_{j=k+1}^{\infty}\widehat{\nu_{2j-1}}(\xi) .
\end{equation}
Here we recall that $\nu_{2j-1} =
 \frac{1}{N_{2j-1}}\sum_{r=0}^{N_{2j-1}-1}\delta_{\frac{r}{N_1\cdots N_{2j-1}}}$
and its Fourier transform is given by
\begin{equation}\label{eq4.4+}
\widehat{\nu_{2j-1}}(\xi) = e^{-\pi i (N_{2j-1}-1)\xi/(N_1\cdots N_{2j-1})}\frac{\sin(\pi \xi/(N_1\cdots N_{2j-2}))}{N_{2j-1}\sin({\pi \xi}/{(N_1\cdots N_{2j-1})})}.
\end{equation}
Let
 $$
 \Lambda_{\mu} = \bigoplus_{j=1}^{\infty}A_{2j-1} = \bigcup_{k=1}^{\infty}\Lambda_k
 $$
(Only finite sums of elements of $A_{2j-1}$, $j\geq 1$, appear in $\Lambda_{\mu}$).
The exponentials $\{e^{2\pi i \lambda x}\}_{\lambda\in\Lambda_{\mu}}$ are mutually
orthogonal in $L^2(\mu)$ by Proposition \ref{prop4.1}. Our goal is verify (\ref{eq4.0}). To do this, we note that, as $Q$ is an entire function, we just need to show that $Q(\xi)\equiv 1$ on a neighborhood of $0$. Let
$$
Q_{k}(\xi) = \sum_{\lambda\in\Lambda_{k}}|\widehat{\mu}(\xi+\lambda)|^2.
$$
 Now, we fix two positive integers $n$ and $p$. By (\ref{eq4.4}) and the fact that
$\{\Lambda_{k}\}_{k\ge 1}$ is an increasing sequence of sets,
\begin{equation}\label{Q}
\begin{aligned}
Q_{n+p}(\xi) =& Q_{n}(\xi)+\sum_{\lambda\in\Lambda_{n+p}\setminus
\Lambda_{n}}|\widehat{\mu}(\xi+\lambda)|^2\\
=&Q_{n}(\xi)+\sum_{\lambda\in\Lambda_{n+p}\setminus
\Lambda_{n}}|\widehat{\mu^{(n+p)}}(\xi+\lambda)|^2\cdot
\left|\prod_{j=n+p+1}^{\infty}\widehat{\nu_{2j-1}}(\xi+\lambda)\right|^2.\\
\end{aligned}
\end{equation}
We need the following proposition which provides a
crucial estimate for the last term in the previous expression
in order to establish the spectral property.

\begin{Prop}\label{prop4.2}
There exists $c>0$ such that
$$
\inf_{k\geq1}\inf_{\lambda\in\Lambda_{k}}\left|\prod_{j=k+1}^{\infty}\widehat{\nu_{2j-1}}(\xi+\lambda)\right|^2\geq c
$$
for all $|\xi|<1/2$, where $\Lambda_{k}$ is given in (\ref{eq4.1}).
\end{Prop}

\begin{pf}
Let $\lambda\in\Lambda_{k}$ and $x_{k,\lambda} = \frac{\xi+\lambda}{N_1N_2\cdots N_{k}}$. We first note that, by (\ref{eq4.4+}),

\begin{equation}\label{eq4.5}
\begin{aligned}
\left|\prod_{j=k+1}^{\infty}\widehat{\nu_{2j-1}}(\xi+\lambda)\right|^2 
 =&\prod_{j=k+1}^{\infty}\frac{\sin^2(\pi ({\xi+\lambda})/({N_1\cdots N_{2j-2}}))}{N_{2j-1}^2\sin^2((\pi (\xi+\lambda))/({N_1\cdots N_{2j-1}}))}\\
 =&\prod_{j=k+1}^{\infty}\frac{\sin^2(\pi x_{2j-2,\lambda})}{N_{2j-1}^2\sin^2(\pi x_{2j-1,\lambda})}.
\end{aligned}
\end{equation}

\medskip

Writing $\lambda = \sum_{j=1}^{k}r_jN_1N_2...N_{2j-2}$ with $0\leq r_j\leq N_{2j-1}-1$, we see immediately that $\lambda\leq N_1\cdots N_{2k-1}-1$. Hence, we have
$$
\frac{\lambda}{N_1\cdots N_{2k}}\leq
\frac{N_1\cdots N_{2k-1}-1}{N_1\cdots N_{2k}}\leq \frac{1}{N_{2k}}\leq \frac12.
$$
Therefore, for all $|\xi|<1/2$, we have
$$
C:=\sup_{k\geq 1}\sup_{\lambda\in\Lambda_{k}}x_{2k,\lambda}=\sup_{k\geq 1}\sup_{\lambda\in\Lambda_{k}}\frac{\xi+\lambda}{N_1\cdots N_{2k}}<\frac{3}{4}
$$
as all $N_j\geq 2$. Note that $N_kx_{k,\lambda} = x_{k-1,\lambda}$ and using two elementary inequalities $\sin x\leq x$ and $\sin x\geq x-\frac{x^3}{3!}$, we have the following estimation for the product in (\ref{eq4.5}),
$$
\begin{aligned}
\prod_{j=k+1}^{\infty}\frac{\sin^2(\pi x_{2j-2,\lambda})}{N_{2j-1}^2\sin^2(\pi x_{2j-1,\lambda})}
\geq&\prod_{j=k+1}^{\infty}\left(1-\frac{\pi^2}{6}x_{2j-2,\lambda}^2\right)^2\\
=&\prod_{j=k+1}^{\infty}\left(1-\frac{\pi^2}{6}\left(\frac{x_{2k,\lambda}}{N_{2k+1}...N_{2j-2}}\right)^2\right)^2\\
\geq& \prod_{j=k+1}^{\infty}\left(1-\frac{\pi^2}{6}\left(\frac{C}{2^{2(j-k)-2}}\right)^2\right)^2\\
=&\prod_{j=1}^{\infty}\left(1-\frac{3\pi^2}{32}\left(\frac{1}{2^{2j-2}}\right)^2\right)^2:=c.
\end{aligned}
$$
As $\sum_{j=1}^{\infty}1/2^{2j-2}<\infty$ and all factors are positive,  $c>0$  and hence the proof is complete.
\end{pf}

\medskip

\noindent{\it Proof of Theorem \ref{th0.1} on ${\Bbb R}^1$.} In view of Theorem \ref{th1.1}, we just need to show that all natural complementary pairs are spectral measures. Let ${\mathcal N}$ be a sequence of positive integers  greater than or equal to 2. If the pair is of Type I, then Proposition \ref{prop4.1} shows that both factors are spectral measures.

\medskip

It remains to consider the Type II case. Let $\mu_{\mathcal N}$ and $\nu_{\mathcal N}$ be defined in (\ref{eq0.2}) and (\ref{eq0.3}). As mentioned before, we only need to prove that  $\mu = \mu_{\mathcal N}$ is a spectral measure.  Let $c$ be the positive number determined in Proposition \ref{prop4.2}.
By Proposition \ref{prop4.1} and (\ref{eq4.0}), we have
$$
\sum_{\lambda\in\Lambda_{n+p}\setminus\Lambda_{n}}|\widehat{\mu^{(n+p)}}(\xi+\lambda)|^2 = 1-\sum_{\lambda\in\Lambda_{n}}|\widehat{\mu^{(n+p)}}(\xi+\lambda)|^2.
$$
Using this fact and  Proposition \ref{prop4.2},  we obtain from (\ref{Q}) that
$$
Q_{n+p}(\xi)\geq Q_{n}(\xi)+c\cdot\left(1-\sum_{\lambda\in\Lambda_{n}}|\widehat{\mu^{(n+p)}}(\xi+\lambda)|^2\right).
$$
Fixing $n$ and letting $p$ go to infinity, it follows that
$$
Q(\xi) \geq Q_{n}(\xi)+c (1-\sum_{\lambda\in\Lambda_{n}}|\widehat{\mu}(\xi+\lambda)|^2) = Q_{n}(\xi)+c(1-Q_{n}(\xi)).
$$
Finally, taking $n$ to infinity, we  obtain that $c(1-Q(\xi))\leq 0$. But $c>0$ and $Q(\xi)\leq 1$
 because $\{e^{2\pi i \lambda x}\}_{\lambda\in\Lambda}$ is an orthogonal set in $L^2(\mu)$.
This show that $Q(\xi)=1$ for $|\xi|\le 1/2$ and thus for all $\xi\in \mathbb{R}$ by analyticity,
completing the proof.

\bigskip

We now establish the tiling property of the spectra. Suppose that we are given a type I decomposition.
Then Proposition \ref{prop4.1} implies that   $\mu_{\mathcal N}$ and $\nu_{\mathcal N}$ have the following spectra:
$$
\Lambda_{\mu} = \bigoplus_{j=1}^{k} A_{2j-1}, \
\Lambda_{\nu} = \bigoplus_{j=1}^{k-1} A_{2j}\oplus N_1\cdots N_{2k-1}{\Bbb Z}.
$$
 It can be seen immediately that $\Lambda_{\mu}\oplus\Lambda_{\nu}= \{0,1,\cdots, N_{2k-1}-1\}\oplus N_{2k-1}{\Bbb Z} = {\Bbb Z}.$

\medskip

 Suppose now the decomposition is of type II.  Note that the complementary measures have the following spectra using the above notations.
$$
\Lambda_{\mu} = \bigoplus_{j=1}^{\infty} A_{2j-1}, \
\Lambda_{\nu} = \bigoplus_{j=1}^{\infty} A_{2j}
$$
Note that $-\Lambda_{\nu}$ is also spectrum of $\nu$. We now claim that $\Lambda_{\mu}\oplus(-\Lambda_{\nu}) = {\Bbb Z}.$  Observe that
$$
A_1\oplus(-A_2) = \{-N_1N_2+N_1,.., N_1-1\}.
$$
$$
A_1\oplus(-A_2)\oplus A_3= \{-N_1N_2+N_1,.., N_1N_2N_3-N_1N_2+N_1-1\}.
$$
Inductively, the sets $A_1\oplus(-A_2)\oplus...\oplus (-1)^{k-1}A_{k}$ cover an
 increasing sequence of consecutive integers.
showing that $\Lambda_{\mu}\oplus(-\Lambda_{\nu}) = {\Bbb Z}$. This proves our claim.
\eproof

\bigskip

\section{Generalized Fuglede's conjecture}

In this section, we will formulate a generalization of Fuglede's
conjecture and prove that it implies the original one. Recall the conjecture we are interested in:

\medskip

\noindent{\bf Conjecture (Generalized Fuglede's Conjecture):}  A compactly supported Borel probability measure $\mu$ on ${\Bbb R}^1$ is spectral if and only if there exists a Borel probability measure $\nu$ and a fundamental domain $Q$ of some lattice on ${\Bbb R}^1$ such that $\mu\ast\nu ={\mathcal L}_{Q}$.

\bigskip

We first prove the following proposition.

\begin{Prop}\label{prop5.1}
Let $\Omega$ and $Q$ be bounded measurable sets of positive Lebesgue measure on ${\Bbb R}^1$.
Suppose that ${\mathcal L}_{\Omega}\ast\nu = {\mathcal L}_Q$, for some Borel probability measure $\nu$. Then
$$
\nu = \sum_{k=1}^{N}\frac{1}{N}\delta_{a_k}, \ Q = \bigcup_{k=1}^{N}(\Omega+a_k)
$$
and ${\mathcal L}((\Omega+a_k)\cap(\Omega+{a_\ell{}})) =0$ for all $k\neq\ell$.
\end{Prop}

\medskip

\begin{pf}
%
%
We first note that ${\mathcal L}_{\Omega}\ast\nu = {\mathcal L}_Q$ if and only if
 $({\mathcal L}_{\Omega}\ast\delta_y)\ast(\nu\ast\delta_x\ast\delta_{-y}) =
 ({\mathcal L}_Q\ast \delta_x)$ for any real numbers $x$ and $y$. Therefore,
there is no loss of generality to assume that the smallest closed intervals containing
$\Omega$ and $Q$ are respectively  $[0,a]$ and $[0,b]$. As $\overline{Q}=$ supp
 $({\mathcal L}_{\Omega}\ast\nu)$ = $\overline{\Omega}$ + supp $\nu$, The support of $\nu$
has to be contained in the non-negative part of the real line.

\medskip

Let $\epsilon>0$ and consider the interval
$E_{\epsilon}=[0,\epsilon)$. Let $\eta_{\epsilon}\in E_{\epsilon}$
be a Lebesgue point of $\chi_{Q}$. Then, using ${\mathcal L}_{\Omega}\ast\nu = {\mathcal L}_Q$,
$$
\begin{aligned}
\frac{1}{{\mathcal L} (Q)}{\mathcal L}\left( Q\cap [\eta_{\epsilon},
\eta_{\epsilon}+h)\right) &= \frac{1}{{\mathcal L}(\Omega)}
\int_{0}^{\eta_{\epsilon}+h} {\mathcal L}\left(\Omega\cap ([\eta_{\epsilon},\eta_{\epsilon}+h)-y)\right)
d\nu(y)\nonumber\\
&=\frac{1}{{\mathcal L}(\Omega)}\int_{0}^{\eta_{\epsilon}+h} {\mathcal L}
\left((\Omega+y)\cap [\eta_{\epsilon},\eta_{\epsilon}+h)\right)d\nu(y),
\end{aligned}
$$
since  $\Omega$ and supp $\nu$ are contained in $ [0,\infty)$. This implies that
$$
\frac{{\mathcal L}(\Omega)}{{\mathcal L}(Q)}{\mathcal L}
(Q\cap [\eta_{\epsilon},\eta_{\epsilon}+h))\leq {\mathcal L}([\eta_{\epsilon},
\eta_{\epsilon}+h))\nu ([0,\eta_{\epsilon}+h))=h\nu ([0,\eta_{\epsilon}+h)).
$$
Since $\eta_{\epsilon}$ is a Lebesgue point of $\chi_{Q}$, we have
$\lim_{h\rightarrow0}\frac{{\mathcal L}(Q\cap[\eta_{\epsilon},\eta_{\epsilon}+h))}{h}=1$.
Therefore, by taking $h\rightarrow 0$, we deduce that
$\frac{{\mathcal L}(\Omega)}{{\mathcal L}(Q)}\leq \nu ([0,\eta_{\epsilon}]).$
Letting $\epsilon$ approach zero, we obtain the inequality
\begin{equation}\label{eq5.2}
\frac{{\mathcal L}(\Omega)}{{\mathcal L}(Q)}\leq \nu(\{0\}).
\end{equation}
Since ${\mathcal L}(\Omega)>0$, $\nu$ has an atom at $0$ and we can write
\begin{equation}\label{eq5.2+}
\nu = p_0\delta_0+(1-p_0)\nu_1, \ p_0 = \nu(\{0\}) \ \mbox{and} \ \nu_1(\{0\})=0.
\end{equation}
The equation ${\mathcal L}_{\Omega}\ast\nu = {\mathcal L}_{Q}$ can thus be rewritten as
\begin{equation}\label{eq5.3}
(1-p_0){\mathcal L}_{\Omega}\ast\nu_1 = {\mathcal L}_{Q}-p_0{\mathcal L}_{\Omega}.
\end{equation}
Since the left hand side of (\ref{eq5.3}) is still a positive measure, this implies that
$$
0\leq ({\mathcal L}_{Q}-p_0{\mathcal L}_{\Omega})(\Omega) \leq \frac{{\mathcal L}(\Omega)}{{\mathcal L}(Q)}-p_0.
$$
Combining it with (\ref{eq5.2}), we conclude that
$p_0=\frac{{\mathcal L}(\Omega)}{{\mathcal L}(Q)}$ and, using (\ref{eq5.3}), we obtain
$$
{\mathcal L}_{\Omega}\ast\nu_1 = {\mathcal L}_{Q\setminus\Omega}
$$

\medskip

If $p_0=1$, then $Q=\Omega$ and $\nu = \delta_0$, so we are done. If not,  we then repeat the argument with
 $Q$ replaced by $Q\setminus\Omega$. We can find $\Omega+a_1\subset Q\setminus\Omega$ such
that $p_1:=\nu_1(\{a_1\})>0$ and $\nu_1 = p_1\delta_{a_1}+(1-p_1)\nu_2$.
Moreover, $p_1 = {\mathcal L}(\Omega)/{\mathcal L}(Q\setminus\Omega)$.
By (\ref{eq5.2+}),
$$
\nu = \frac{{\mathcal L} (\Omega)}{{\mathcal L} (Q)}\left(\delta_0+\delta_{a_1}\right)+(1-p_1)\nu_2.
$$
The theorem will be proved if $p_1=1$. Otherwise, we continue this process to
obtain a maximal number $N$ of measure disjoint translates of $\Omega$,
$\Omega+a_1$,..,$\Omega+a_{N-1}$ such
that $Q\supset\bigcup_{k=0}^{N-1}(\Omega+a_k)$. Since ${\mathcal L}(\Omega)>0$ and
${\mathcal L}(Q)\geq N{\mathcal L}(\Omega)$,  $N$ is the largest integer such
that ${\mathcal L}(Q)\geq N{\mathcal L}(\Omega)$.  We can then write
$$
\nu = \frac{{\mathcal L}(\Omega)}{{\mathcal L}(Q)}\left(\delta_0+...+\delta_{a_{N-1}}\right)+(1-p_{N-1})\nu_{N}.
$$
If $p_{N-1}<1$, we could iterate this process to obtain one more disjoint translate of $\Omega$ contained in $Q$, which is certainly impossible by this choice of $N$. Hence, $p_{N-1}=1$. As $\nu$ is a probability measure, we must have ${\mathcal L}(\Omega)/{\mathcal L}(Q) = 1/N$. Therefore, the proposition is proved.
\end{pf}

\medskip

\begin{theorem}
The validity of generalized Fuglede's conjecture implies that of the original Fuglede's conjecture on ${\Bbb R}^1$.
\end{theorem}

\begin{pf}
Suppose that $\Omega$ is a bounded spectral set, then ${\mathcal L}_{\Omega}$ is a spectral measure. By the generalized
 Fuglede's conjecture, we can find a probability measure $\nu$ and a fundamental domain $Q$ of some lattice $\Gamma$ such that
$$
{\mathcal L}_{\Omega}\ast\nu = {\mathcal L}_{Q}.
$$
By Proposition \ref{prop5.1}, $\nu$ is a purely discrete measure that can be written as $\nu = \frac{1}{\#{\mathcal A}}\delta_{\mathcal A}$ for some finite discrete subset ${\mathcal A}$ and
$$
Q = \bigcup_{a\in{\mathcal A}}(\Omega+a).
$$
As $Q$ is a fundamental domain $Q$ of the lattice $\Gamma$, $\Omega$ is a translational tile with tiling set given by ${\mathcal A}+\Gamma$.

\medskip

Conversely, suppose that $\Omega$ is a bounded translational tile with tiling set ${\mathcal J}$. By the result of Lagarias and Wang \cite{[LW1]}, all tiling sets on ${\Bbb R}^1$ are periodic. This implies that we can find a finite set $A\subset{\Bbb R}$ and a lattice $\Gamma$ such that ${\mathcal J} = {\mathcal A}+\Gamma$. This means that  the set $Q = \Omega+{\mathcal A}$ is a fundamental domain of $\Gamma$. Letting $\nu = \frac{1}{\#{\mathcal A}}\delta_{\mathcal A}$, ${\mathcal L}_{\Omega}\ast\nu = {\mathcal L}_{Q}$. By the
generalized Fuglede's conjecture, ${\mathcal L}_{\Omega}$ is a spectral measure and $\Omega$ is a spectral set.
\end{pf}

\bigskip

\section{The Higher Dimensional Case}

Let $\mu_1$,...,$\mu_d$ be Borel probability measures on ${\Bbb R}^1$. The Cartesian product of these measures is the unique Borel probability measure $ \mu_{1}\otimes...\otimes\mu_d$ on ${\Bbb R}^d$ such that
$$
(\mu_{1}\otimes...\otimes\mu_d)(E_1\times...\times E_d) = \prod_{i=1}^{d}\mu_i(E_i),
$$
for any Borel sets $E_i$, $1\leq i\leq d$, on ${\Bbb R}^1$. In this section, we characterize the measures $\mu$ and $\nu$ on ${\Bbb R}^d$ which are solutions of the equation
\begin{equation}\label{eq2.1}
\mu\ast\nu = {\mathcal L}_{[0,1]^d}.
\end{equation}
as Cartesian products of the measures satisfying the corresponding one-dimensional equation.

\begin{theorem}\label{th2.2}
Let $\mu$ and $\nu$ be compactly supported probability measures on ${\Bbb R}^d$. Then  $\mu$ and $\nu$ are solutions to (\ref{eq2.1}) if and only if there exists compactly supported Borel probability measures $\{\sigma_i\}_{i=1}^{d}$ and $\{\tau_i\}_{i=1}^{d}$ on ${\Bbb R}^1$ such that
\begin{equation}\label{eq2.3}
\mu = \sigma_{1}\otimes...\otimes\sigma_d, \ \nu = \tau_{1}\otimes...\otimes\tau_d
\end{equation}
and $\sigma_i\ast\tau_i = {\mathcal L}_{[0,1]}$ for all $i=1,...,d$.
\end{theorem}

Note that the sufficiency part of the theorem follows by a direct computation.
We only need to establish the necessity part of the theorem.
 Denote by $P$ the orthogonal projection of the first coordinate on ${\Bbb R}^d$ and $Q$ the orthogonal projection of the corresponding orthogonal complement. If $\mu$ is a positive Borel measure on ${\Bbb R}^d$, we denote by $\mu P^{-1}$ the positive Borel measure on ${\Bbb R}^1$ defined by
 $\mu P^{-1}(E)= \mu(P^{-1}(E))$ for any Borel set $E\subset{\Bbb R}$ and the measure $\mu Q^{-1}$ is similarly defined.  We will need the following lemmas.

\medskip

\begin{Lem}\label{Lem2.1}
Let $\mu$ and $\nu$ be two probability measures on ${\Bbb R}^d$. Then
$$
(\mu\ast\nu) P^{-1} = (\mu P^{-1})\ast(\nu P^{-1}), \ \mbox{and} \ (\mu\ast\nu) Q^{-1} = (\mu Q^{-1})\ast(\nu Q^{-1}).
$$
In particular, if $\mu$ and $\nu$ are two Borel probability measures satisfying (\ref{eq2.1}), then we have
$$
(\mu P^{-1})\ast(\nu P^{-1}) = {\mathcal L}_{[0,1)} \ \mbox{and} \ (\mu Q^{-1})\ast(\nu Q^{-1}) = {\mathcal L}_{[0,1]^{d-1}}.
$$
\end{Lem}

\medskip

\begin{pf} The proof follows easily from the fact that
$$
(\mu P^{-1})^{\widehat{}}(\xi) = \widehat{\mu}(\xi,0,...,0), \ \mbox{and} \ (\mu Q^{-1})^{\widehat{}}(\xi_2,...,\xi_d) = \widehat{\mu}(0,\xi_2,...,\xi_d).
$$
\end{pf}
%

\medskip

\begin{Lem}\label{Lem2.4}
Let $\nu$ be a Borel probability measure on ${\Bbb R}^d$. Then, there is at most one probability measure $\mu$ on ${\Bbb R}^d$ satisfying $\mu\ast\nu = {\mathcal L}_{[0,1]^d}$.
\end{Lem}

\medskip

\begin{pf}
If $\mu$ is as above, we have
\begin{equation}\label{eq5.0}
\widehat{\mu}(\xi)\widehat{\nu}(\xi)=\left({\mathcal L}_{[0,1]^d}\right)^{\widehat{}}(\xi), \ \xi\in{\Bbb R}^d.
\end{equation}
Therefore, $\widehat{\mu} (\xi)$ is thus determined on the set
$$
F = \{\xi\in{\Bbb R}^d: \xi_i\not\in{\Bbb Z}^d, i=1,...,d\},
$$
Since $\overline{F}={\Bbb R}^d$ and $\widehat{\mu}$ is continuous (as $\mu$ is compactly supported), $\widehat{\mu}$ and thus $\mu$ is completely determined by (\ref{eq5.0}).
\end{pf}

\medskip

The previous lemma is also valid if $[0,1]^d$ is replaced by a $d$-dimensional rectangular box.
Now, we proceed to the proof of Theorem \ref{th2.2}.

\medskip

\noindent{\it Proof of Theorem \ref{th2.2}.} We prove the necessity part of the theorem by induction on the dimension. The statement is proved when $d=1$ in Theorem \ref{th1.1}. Assuming that the statement is true for $d-1$, we now establish  it  on ${\Bbb R}^d$.

\medskip

Let $\mu$ and $\nu$ be two Borel probability measures satisfying $\mu\ast\nu = {\mathcal L}_{[0,1]^d}$. By Lemma \ref{Lem2.1} and Theorem \ref{th1.1} (see also equation (\ref{eq1.2})), we can find an integer $N_1\geq 2$ such that $\mu P^{-1}$ and $\nu P^{-1}$ can be decomposed (after possibly interchanging these two measures) as
\begin{equation}\label{eq5.5}
\mu P^{-1} = \nu_1\ast\alpha_1, \ \mbox{and} \  \alpha_1\ast(\nu P^{-1}) = {\mathcal L}_{[0,1/N_1]}
\end{equation}
where $\nu_1 = 1/N_1\sum_{j=0}^{N_1-1}\delta_{j/N_1}$ and $\alpha_1(E) = N_1 (\mu P^{-1})(E\cap[0,1/N_1))$ for any Borel set $E$. Let $C_{N_1}$ be the $d$-dimensional rectangular box $\left[0,\frac{1}{N_1}\right)\times[0,1]^{d-1}$. Then $[0,1]^d\setminus C_{N_1} = \left[\frac{1}{N_1},1\right]\times[0,1]^{d-1}$ and
$$
\mu\left(C_{N_1}\right) = \mu P^{-1}\left(\left[0,\frac{1}{N_1}\right)\right) = \frac{1}{N_1}.
$$
Hence, we can define  two Borel probability measures  on ${\Bbb R}^d$, $\rho_1$ and $\widetilde{\rho_1}$, satisfying
$$
\rho_1(E) =N_1 \mu\left(E\cap C_{N_1}\right), \ \widetilde{\rho_1}(E) = \frac{N_1}{N_1-1}\mu\left(E\cap \left([0,1]^d\setminus C_{N_1}\right)\right)
$$
for any Borel sets $E$. Then $\mu =\frac{1}{N_1}\rho_1+(1-\frac{1}{N_1})\widetilde{\rho_1}$. Since supp $\widetilde{\rho}\subset [0,1]^d\setminus C_{N_1}$ and supp $\nu\subset [0,1]^d$, we have $\nu\ast\widetilde{\rho}=0$ on the rectangular box $C_{N_1}$. Hence,
$$
\rho_1\ast\nu =N_1(\mu\ast\nu) = {\mathcal L}_{C_{N_1}} \  \ \mbox{on}   \ \ C_{N_1}.
$$
We can thus write $\rho_1\ast\nu = {\mathcal L}_{C_{N_1}}+\eta$ where $\eta$ is a positive measure. However, $\eta=0$ as $\rho_1\ast\nu$ and ${\mathcal L}_{C_{N_1}}$ are probability measures. Hence,
$$
(\nu_1\otimes\delta_{0_{d-1}})\ast\rho_1\ast\nu =\left(\frac{1}{N_1}\sum_{j=0}^{N_1-1}\delta_{(j/N_1,0...,0)}\right)\ast\rho_1\ast\nu = {\mathcal L}_{[0,1]^d}
$$
where $0_{d-1} = (0,...,0)\in{\Bbb R}^{d-1}$. By Lemma \ref{Lem2.4}, we have that
\begin{equation}\label{eq5.6}
\mu = (\nu_1\otimes\delta_{0_{d-1}})\ast\rho_1, \ \mbox{and} \ \rho_1\ast\nu = {\mathcal L}_{[0,1/N_1]\times[0,1]^{d-1}}
\end{equation}
Furthermore, $\rho_1 P^{-1} = \alpha_1 $ where $\alpha_1$ is defined in (\ref{eq5.5}).

\medskip

We now consider two cases depending on whether $\mu P^{-1}$ and $\nu P^{-1}$ correspond to a type I or type II decomposition (as defined in Section 2).

\medskip

\noindent{\bf Case 1 (Type I decomposition):} Using the notations introduced in Section 2, we have then, without loss of generality, that
 $$
 \mu P^{-1} = \nu_1\ast...\nu_{2k-1}, \ \nu P^{-1} = \nu_2\ast...\nu_{2k}\ast{\mathcal L}_{[0,\frac{1}{N_1...N_{2k}}]}.
 $$
 By the previous steps, the identities in (\ref{eq5.6}) hold. A similar argument, shows  the existence of a probability measure $\rho_2$  such that
 $$
 \nu = (\nu_2\otimes\delta_{0_{d-1}})\ast\rho_2 \ \  \mbox{and} \ \ \rho_1\ast\rho_2 = {\mathcal L}_{[0,\frac{1}{N_1N_2}]\times[0,1]^{d-1}}.
 $$
Continuing this procedure $2k$-times, we deduce the existence of probability measures $\rho_{2k-1}$ and $\rho_{2k}$ such that
\begin{equation}\label{eq5.8}
\mu = ((\nu_1\ast\nu_3\ast...\ast\nu_{2k-1})\otimes\delta_{0_{d-1}})\ast\rho_{2k-1}
\end{equation}
\begin{equation}\label{eq5.9}
\nu  = ((\nu_2\ast\nu_4\ast...\ast\nu_{2k})\otimes\delta_{0_{d-1}})\ast\rho_{2k}
\end{equation}
and
\begin{equation}\label{eq5.10}
\rho_{2k-1}\ast\rho_{2k} = {\mathcal L}_{[0,1/N_1N_2...N_{2k}]\times[0,1]^{d-1}}.
\end{equation}
By (\ref{eq5.8}) and Lemma \ref{Lem2.1}, $\mu P^{-1} = \nu_1\ast...\nu_{2k-1}\ast\rho_{2k-1} P^{-1}$, showing that $\rho_{2k-1} P^{-1} =\delta_0$. Hence, we can write $\rho_{2k-1} = \delta_{0}\otimes\sigma$ for some positive measure $\sigma$ on ${\Bbb R}^{d-1}$. Using (\ref{eq5.10}) and Lemma \ref{Lem2.1} again, we obtain that $\sigma\ast(\rho_{2k} Q^{-1}) = {\mathcal L}_{[0,1]^{d-1}}$. Hence,
$$
\begin{aligned}
\rho_{2k-1}\ast\rho_{2k} =& {\mathcal L}_{[0,1/N_1N_2...N_{2k}]}\otimes{\mathcal L}_{[0,1]^{d-1}} \\
=& {\mathcal L}_{[0,1/N_1N_2...N_{2k}]}\otimes(\sigma\ast(\rho_{2k} Q^{-1}) )\\
 =& (\delta_0\otimes \sigma)\ast({\mathcal L}_{[0,1/N_1N_2...N_{2k}]}\otimes(\rho_{2k} Q^{-1}))\\
 =&\rho_{2k-1}\ast({\mathcal L}_{[0,1/N_1N_2...N_{2k}]}\otimes(\rho_{2k} Q^{-1})).
\end{aligned}
$$
 Lemma \ref{Lem2.4} shows that $\rho_{2k} = {\mathcal L}_{[0,1/N_1N_2...N_{2k}]}\otimes(\rho_{2k} Q^{-1})$ and (\ref{eq5.9}) implies that $
\nu = \nu P^{-1}\otimes \rho_{2k} {Q}^{-1}.$  Finally, applying the induction hypothesis to the identity $\sigma\ast(\rho_{2k} Q^{-1}) = {\mathcal L}_{[0,1]^{d-1}}$, we can write $\sigma = \sigma_2\otimes...\otimes\sigma_d$ and $\rho_{2k-1} Q^{-1} = \tau_2\otimes...\otimes\tau_d$ with $\sigma_i\ast\tau_i = {\mathcal L}_{[0,1]}$ and Theorem \ref{th2.2} for dimension $d$ follows.

\medskip

%

\noindent{\bf Case 2 (Type II decomposition).} In this case, we can without loss of generality assume that
$$
\mu P^{-1} = \nu_1\ast\nu_3\ast..., \ \nu P^{-1} = \nu_2\ast\nu_4\ast...
$$
 and we still have (\ref{eq5.8}), (\ref{eq5.9}) and (\ref{eq5.10}) for all $k= 1,2,....$
with $\rho_n P^{-1} \ne \delta_0$ for any integer $n$. As $\rho_{n}$ are all probability measures,
we can assume,
by passing to subsequences if necessary, that the sequences $\{\rho_{2k-1}\}$ and $\{\rho_{2k}\}$
converge weakly to some probability measures that we denote by $\sigma$ and $\tau$, respectively.
From (\ref{eq5.10}), it is immediate to see that the supports of $\sigma$ and $\tau$ are both contained in $\{0\}\times[0,1]^{d-1}$. We can write $\sigma = \delta_0\otimes\sigma'$ and $\tau = \delta_0\otimes \tau'$.  By passing to weak limit in  (\ref{eq5.8}) and (\ref{eq5.9}), we have
\begin{equation}\label{eq5.11}
\mu = (\mu P^{-1}\otimes\delta_{0_{d-1}})\ast (\delta_0\otimes\sigma'), \  \ \nu = (\nu P^{-1}\otimes\delta_{0_{d-1}})\ast (\delta_0\otimes\tau').
\end{equation}
As $\mu\ast\nu = {\mathcal L}_{[0,1]^d}$ and $(\mu P_1^{-1}\otimes\delta_{0_{d-1}})\ast  (\nu P_1^{-1}\otimes\delta_{0_{d-1}}) = {\mathcal L}_{[0,1]}\otimes\delta_{0_{d-1}}$, we have
$$
\sigma'\ast\tau '={\mathcal L}_{[0,1]^{d-1}},
$$
The conclusion follows immediately by (\ref{eq5.11}) using the induction hypothesis.
\eproof
\medskip

\noindent{\it Proof of Theorem \ref{th0.1} on ${\Bbb R}^d$.} The proof follows from the result on ${\Bbb R}^1$. By Theorem \ref{th2.2}, we can write $\mu = \sigma_1\otimes...\otimes\sigma_d$ and $\nu = \tau_1\otimes...\otimes\tau_d$ with $\sigma_i\ast\tau_i = {\mathcal L}_{[0,1]}$. Therefore, our conclusion on ${\Bbb R}^1$ implies that $\sigma_i$ and $\tau_i$ are spectral measures on ${\Bbb R}^1$ with spectrum $\Lambda_{\sigma_i}$ and $\Lambda_{\tau_i}$ respectively. Moreover, they satisfies $\Lambda_{\sigma_i}\oplus\Lambda_{\tau_i} = {\Bbb Z}$. Now  we define
$$
\Lambda_{\mu} = \bigotimes_{i=1}^{d}\Lambda_{\sigma_i}, \ \Lambda_{\nu} = \bigotimes_{i=1}^{d}\Lambda_{\tau_i},
$$
where $\bigotimes_{i=1}^{d} A_i := \{(a_1,...,a_d): a_i\in A_i\}$ for sets $A_i\subset{\Bbb R}^1$.
We claim that $\Lambda_{\mu}$ is a spectrum for $\mu$ (the proof that $\Lambda_{\nu}$ is a spectrum for $\nu$ is similar).

\medskip

Note that $
\widehat{\mu}(\xi) = \prod_{i=1}^{d}\widehat{\sigma_i}(\xi_i).$
 From this, it follows easily that
$$
\sum_{\lambda\in\Lambda_{\mu}} |\widehat{\mu}(\xi+\lambda)|^2 = \prod_{i=1}^{d}\left(\sum_{\lambda_i\in\Lambda_{\sigma_i}}|\widehat{\sigma_i}(\xi_i+\lambda_i)|^2\right)=1.
$$
Hence, $\Lambda_{\mu}$ is a spectrum for $\mu$. That the tiling property of the spectra (i.e. $\Lambda_{\mu}\oplus\Lambda_{\nu} = {\Bbb Z}^d$) follows immediately from the tiling property of $\Lambda_{\sigma_i}$ and $\Lambda_{\tau_i}$.
\eproof

\bigskip

\section{Remarks and Open questions}

\medskip

As indicated in the introduction, the statement $\bf{{\mathcal F}(Q)}$ is false in general.
Nonetheless, this statement suggests many related questions that may help us understand the
relationship among convolutions, translational tilings and spectral measures.
Motivated by the generalized Fuglede's conjecture,
one of the main questions we would like to ask is:

  \medskip

\noindent {\bf (Q1):} For which $Q$ is the statement ${\bf {\mathcal F}(Q)}$ true?

 \medskip

  This question seems to be hard if we go beyond cubes
 as the methods of this paper would be difficult to extend.
An easier, but still interesting question concerns the decomposition of
the Lebesgue measure on sets as convolution product of singular  measures:

 \medskip

\noindent {\bf (Q2):} For what kind of measurable (resp. spectral) sets $Q$
can ${\mathcal L}_{Q}$ be decomposed
into the convolution of two singularly continuous  (resp. spectral) measure ?

 \medskip

One natural type of such sets will be the self-affine tiles \cite{[LW2]}.
 These tiles can be described as  infinite convolution product of discrete measures and
can therefore be decomposed into two singular measures using  methods similar to those in Section 2.

\medskip

Fourier frames and exponential Riesz bases are natural generalization of exponential orthonormal bases.
It has been an interesting question to produce singular measures with
 Fourier frames but not exponential orthonormal bases.
By now we only know we can produce such measures by considering measures
which are absolutely continuous with respect to a spectral measure
with density bounded above and away from 0  or
convolving a spectral measure with some discrete measures \cite{[HLL], [DL1]}.
These methods are rather restrictive.  As absolutely continuous (w.r.t. Lebesgue) measures
with Fourier frames were completely classified in \cite{[Lai]}, we ask

 \medskip

\noindent {\bf (Q3):} Can we produce new singular measures admitting Fourier
frames by decomposing an absolutely continuous (w.r.t. Lebesgue) measures with Fourier frames?
Conversely, is it true that all measures admitting Fourier frames are constructed in this way?

 \medskip

Given a spectral measure $\mu$, another important issue is to classify its spectrum.
This question has been studied  for Lebesgue measures and some Cantor measures  in
\cite{[LRW],[DHS],[DHL]}. However, there is no satisfactory answer when the measure is
singular. The tiling statement of Theorem \ref{th0.1}, suggests
a possible answer.
 \medskip

\noindent {\bf (Q4):} Let $\mu$ and $\nu$ be a natural complementary pair of
${\mathcal L}_{[0,1]}$. Let also $\Lambda_{\mu}$ be a spectrum for $L^2(\mu)$,
does there exist a spectrum $\Lambda_{\nu}$ for $L^2(\nu)$ such that
$\Lambda_{\mu}\oplus\Lambda_{\nu} = {\Bbb Z}$?

 \medskip

It is not difficult to prove that  {\bf (Q4)} actually holds for type I
decompositions. The remaining challenge is to answer the question for type II decompositions.

\end{document}